\documentclass[11pt]{amsart}

\newtheorem{theorem}{Theorem}
\newtheorem{lemma}[theorem]{Lemma}

\newtheorem{corollary}[theorem]{Corollary}

\theoremstyle{definition}

\newtheorem{example}[theorem]{Example}

\newtheorem{question}[theorem]{Question}

\theoremstyle{remark}

\numberwithin{equation}{section}

\newcommand{\Ker}{\operatorname{ker}}

\begin{document}

\title[A Banach-Stone theorem for Riesz isomorphisms]
{A Banach-Stone theorem for Riesz
isomorphisms of Banach lattices}

\author[J.X. Chen]
{Jin Xi Chen}

\address{Department of Mathematics, Southwest Jiaotong
University, Chengdu 610031, P.R. China}
 \email{jinxichen@home.swjtu.edu.cn}

\author[Z.L. Chen]
{Zi Li Chen}
\address{Department of Mathematics, Southwest Jiaotong
University, Chengdu 610031, P.R. China}
\email{zlchen@home.swjtu.edu.cn}
\author[N.-C. Wong]
{Ngai-Ching Wong}
\address{Department of Applied Mathematics, National Sun Yat-sen University, Kaohsiung 80424, Taiwan}
\email{wong@math.nsysu.edu.tw}


\begin{abstract}
Let $X$ and $Y$ be compact Hausdorff spaces, and $E$, $F$ be
Banach lattices. Let $C(X,E)$ denote the Banach lattice of all
continuous $E$-valued functions on $X$ equipped with the pointwise
ordering and the sup norm. We prove that if there exists a Riesz
isomorphism $\mathnormal{\Phi}: C(X,E)\rightarrow C(Y,F)$ such
that $\mathnormal{\Phi}f$ is non-vanishing on $Y$ if and only if
$f$ is non-vanishing on $X$, then $X$ is homeomorphic to $Y$, and
$E$ is Riesz isomorphic to $F$. In this case, $\mathnormal{\Phi}$
can be written as a weighted composition operator:
$\mathnormal{\Phi} f(y)=\mathnormal{\Pi}(y)(f(\varphi(y)))$, where
$\varphi$ is a homeomorphism from $Y$ onto $X$, and
$\mathnormal{\Pi}(y)$ is a Riesz isomorphism from $E$ onto $F$ for
every $y$ in $Y$. This generalizes some known results obtained
recently.

\end{abstract}

\maketitle

\section{Introduction}
Let $X$ and $Y$ be compact Hausdorff spaces, and $C(X)$, $C(Y)$
denote the spaces of real-valued continuous functions defined on
$X$, $Y$ respectively. There are three versions of the Banach-Stone
theorem. That is to say, surjective linear isometries,  ring
isomorphisms and   lattice isomorphisms from $C(X)$ onto $C(Y)$
yield  homeomorphisms between $X$ and $Y$, respectively (cf. \cite{1, 6, 13}).

\par Jerison \cite{12} got the first vector-valued version of
the Banach-Stone theorem. He proved that if the Banach space $E$ is
strictly convex, then every surjective linear isometry
$\mathnormal{\Phi}:C(X,E)\rightarrow C(Y,E)$ can be written as a
weighted composition operator
 $$
 \mathnormal{\Phi} f(y)=\mathnormal{\Pi } (y)(f(\varphi(y))), \qquad \forall\,f\in C(X, E),\,\forall\,y\in Y.
 $$
 Here $\varphi$ is a homeomorphism from
$Y$ onto $X$, and $\mathnormal{\Pi }$ is a continuous map from $Y$ into the space
 $(\mathcal{L}(E, E), SOT)$ of bounded linear operators on  $E$ equipped with
 the strong operator topology (\textit{SOT}).
Furthermore,  $\mathnormal{\Pi}(y)$ is a surjective linear isometry on $E$ for every $y$ in $Y$.
After Jerison \cite{12}, many vector-valued versions of the Banach-Stone theorem
have been obtained in different ways (see, e.g., \cite{3, 4, 5, 7, 9, 10,
11, 14}).

Let $E$, $F$ be  nonzero real Banach lattices, and $C(X,E)$ be the
Banach lattice of all continuous $E$-valued functions on $X$
equipped with the pointwise ordering and the sup norm. Note that,
in general, a Riesz isomorphism (\,i.e., lattice isomorphism\,)
from $C(X,E)$ onto $C(Y,F)$ does not necessarily induce a
topological homeomorphism from $X$ onto $Y$ (cf. \cite[Example
3.5]{14}). To consider the Banach-Stone theorems for continuous
Banach lattice valued functions, we would like to mention the
papers \cite{5, 7, 14}. In particular, when $E$, $F$ are both
Banach lattices and Riesz algebras, Miao, Cao and Xiong \cite{14}
recently proved that if $F$ has no zero-divisor and there exists a
Riesz algebraic isomorphism $\mathnormal{\Phi}: C(X,E)\rightarrow
C(Y,F)$ such that $\mathnormal{\Phi}f$ is non-vanishing on $Y$ if
$f$ is non-vanishing on $X$, then $X$ is homeomorphic to $Y$, and
$E$ is Riesz algebraically isomorphic to $F$. By saying $f$ in
$C(X,E)$ is \emph{non-vanishing},  we mean that $0\notin f(X)$.
Indeed, under these conditions they obtained that
$\mathnormal{\Phi}^{-1}g$ is non-vanishing on $X$ if $g\in C(Y,F)$
is non-vanishing on $Y$. Note that every Riesz algebraic
isomorphism must be a Riesz isomorphism.

\par Let $E$ and $F$ be Banach lattices. More recently, Ercan and \"{O}nal \cite{7} have established that
if $F$ is an \textit{AM}-space with unit, i.e., a $C(K)$-space,
and there exists a Riesz isomorphism $\mathnormal{\Phi}:
C(X,E)\rightarrow C(Y,F)$ such that $\mathnormal{\Phi}f$ is
non-vanishing on $Y$ if and only if $f$ is non-vanishing on $X$,
that is, both $\mathnormal{\Phi}$ and $\mathnormal{\Phi}^{-1}$ are
non-vanishing preserving, then $X$ is homeomorphic to $Y$, and $E$
is Riesz isomorphic to $F$.

\par Inspired by \cite{5, 7, 14}, one can set a natural question:

\begin{question}
               Is $X$ homeomorphic to
$Y$ if $E$, $F$ are Banach lattices and there exists a Riesz
isomorphism $\mathnormal{\Phi}: C(X,E)\rightarrow C(Y,F)$ such that
both $\mathnormal{\Phi}$ and $\mathnormal{\Phi}^{-1}$  are
non-vanishing preserving?
\end{question}

\par In this paper we show the answer to the above question is affirmative. Moreover, in this case $\mathnormal{\Phi}$ can
be written as a weighted composition operator:
$$
\mathnormal{\Phi}
f(y)=\mathnormal{\Pi}(y)(f(\varphi(y))),\qquad \forall\,f\in C(X,
E),\,\forall\,y\in Y,
$$
where $\varphi$ is a homeomorphism from $Y$
onto $X$, and $\mathnormal{\Pi}(y)$ is a Riesz isomorphism from $E$
onto $F$ for every $y$ in $Y$. This generalizes the results obtained
by Cao, Reilly and Xiong \cite{5}, Miao, Cao, and Xiong \cite{14}, and Ercan
and \"{O}nal \cite{7}.
\par Our notions are standard. For the undefined notions and basic facts
concerning Banach lattices we refer the readers to the monographs
\cite{1, 2, 13}.

\section{A Banach-Stone theorem for Riesz isomorphisms}
\par In the following
 we always assume $X$ and $Y$ are compact
Hausdorff spaces, $E$ and $F$ are nonzero Banach lattices, and
$\mathcal{L}(E,F)$ is the space of bounded linear operators from $E$
into $F$ equipped with \textit{SOT}. For $x$ in $X$ and $y$ in $Y$, let
$M_{x}$ and $N_{y}$ be defined as$$M_{x}=\{f\in
C(X,E):\,f(x)=0\},\quad N_{y}=\{g\in C(Y,F): g(y)=0\}.$$Clearly,
$M_{x}$ and $N_{y}$ are closed (order) ideals in $C(X,E)$ and
$C(Y,F)$, respectively.

\begin{lemma}\label{Lemma1}
Let $\mathnormal{\Phi}: C(X,E)\rightarrow C(Y,F)$ be a Riesz
isomorphism  such that $\mathnormal{\Phi}(f)$ is non-vanishing on $Y$ if and
only if $f$ is non-vanishing on $X$. Then for each $x$ in $X$ there exits a
unique $y$ in $Y$ such that
$$
\mathnormal{\Phi} M_{x}=N_{y}.
$$
In particular,
this defines a bijection $\varphi$ from $Y$ onto $X$ by
$\varphi(y)=x$.
\end{lemma}
\begin{proof}
For each $x$ in $X$, let
$$
\mathcal{Z}(\mathnormal{\Phi} M_{x}) = \{y\in Y : \mathnormal{\Phi} f(y)=0\,\,  \text{for all}\, f\in
M_x\}.
$$
We first claim that $\mathcal{Z}(\mathnormal{\Phi} M_{x})$ is
non-empty. Suppose, on the contrary, that
$\mathcal{Z}(\mathnormal{\Phi} M_{x})$ is empty. Then for each $y$
in $Y$ there would exist an $f_{y}$ in $M_{x}$ such that
$\mathnormal{\Phi} f_{y}(y)\neq0$, and thus $\mathnormal{\Phi}
f_{y}$ is non-vanishing in an open neighborhood of $y$. Note that
$|f_{y}|\in M_{x}$, and $\mathnormal{\Phi}
|f_{y}|=|\mathnormal{\Phi} f_{y}|$ since $\mathnormal{\Phi}$ is a
Riesz isomorphism. Therefore, we can assume further that both
$f_{y}$ and $\mathnormal{\Phi} f_{y}$ are positive, by replacing
them by their absolute values if necessary. By the compactness of
$Y$, we can choose finitely many $f_{1}, \ldots, f_{n}$ from
$M_{x}^{+}$ such that the positive functions $\mathnormal{\Phi}
f_{1},\ldots, \mathnormal{\Phi} f_{n}$ have no common zero in $Y$.
Hence $\mathnormal{\Phi}(f_{1}+\cdots + f_n)$ is strictly positive,
that is, $\mathnormal{\Phi}(f_{1}+\cdots + f_{n})(y)>0$ for each $y$
in $Y$. This contradicts the fact that $f_{1}+\cdots + f_{n}$
vanishes at $x$. We thus prove that $\mathcal{Z}(\mathnormal{\Phi}
M_{x})\neq\phi$.

\par Next, we claim that $\mathcal{Z}(\mathnormal{\Phi} M_{x})$ is a singleton. Indeed, if
$y_{1}, y_{2}\in \mathcal{Z}(\mathnormal{\Phi} M_{x})$ then we would
have $\mathnormal{\Phi} M_{x} \subseteq N_{y_{i}}, i=1,2$.
Applying the above argument to $\mathnormal{\Phi}^{-1}$, we shall
have $ \mathnormal{\Phi}^{-1}N_{y_{i}} \subseteq M_{x_{i}}$ for some
$x_{i}$ in $X$, \,$i=1, 2$. It follows that $\mathnormal{\Phi} M_{x}
\subseteq N_{y_{i}}\subseteq \mathnormal{\Phi} M_{x_{i}},\, i=1,2$.
Then $ x = x_{1} = x_{2} $ since  $\mathnormal{\Phi}$ is bijective
and $X$ is Hausdorff. Thus,
$$
y_{1} = y_{2}\quad \text{and }\quad
\mathnormal{\Phi} M_{x} = N_{y_1}=N_{y_2}.
$$
Now, we can define a bijective map $\varphi: Y\rightarrow X$ such
that
$$
 \mathnormal{\Phi} M_{\varphi(y)} = N_{y}, \quad\forall\, y\in Y.$$
\end{proof}

The following main  result answers affirmatively the question
mentioned in the introduction, and solves the conjecture of Ercan
and \"{O}nal in \cite{7}.

\begin{theorem}\label{Theorem2}
Let $\mathnormal{\Phi}: C(X,E)\rightarrow C(Y,F)$ be a Riesz
isomorphism  such that $\mathnormal{\Phi}f $ is non-vanishing on
$Y$ if and only if $f$ is non-vanishing on $X$. Then $Y$ is
homeomorphic to $X$, and $\mathnormal{\Phi}$ can be written as a
weighted composition operator
$$
\mathnormal{\Phi}f(y)=\mathnormal{\Pi}(y)(f(\varphi(y))), \qquad
\forall\,f\in C(X,E),\,\forall\,y\in Y.
$$
Here $\varphi$ is a homeomorphism from $Y$ onto $X$, and
$\mathnormal{\Pi}(y)$ is a Riesz isomorphism from $E$ onto $F$ for
every $y$ in $Y$. Moreover,
$\mathnormal{\Pi}:Y\rightarrow(\mathcal{L}(E, F), SOT)$ is
continuous, and $\|\mathnormal{\Phi}\|=\sup_{y\in
Y}\|\mathnormal{\Pi}(y)\|$.
\end{theorem}
\begin{proof}
First, we show that the bijection $\varphi$ given in Lemma
\ref{Lemma1} is a homeomorphism from $Y$ onto $X$. It suffices to
verify the continuity of $\varphi$ since $Y$ is compact and $X$ is
Hausdorff. To this end, suppose, to the contrary, that there would
exist a net $\{y_{\lambda}\}$ in $Y$ converging to $y_{0}$ in $Y$,
but $\varphi(y_{\lambda})$ converges to $x_{0}\neq\varphi(y_{0})$
in $X$.
\par Let $U_{x_{0}}$ and $U_{\varphi(y_{0})}$ be disjoint open
neighborhoods of $x_{0}$ and $\varphi(y_{0})$, respectively. First,
for any $f$ in $C(X,E)$ vanishing outside $U_{\varphi(y_{0})}$ we
claim that $\mathnormal{\Phi} f(y_{0})=0$. Indeed, since
$\varphi(y_{\lambda})$ belongs to $U_{x_{0}}$ for $\lambda$ large
enough and $f(x)=0$ for any $x$ in $U_{x_{0}}$, we have that $f\in
M_{\varphi(y_{\lambda})}$. It follows from Lemma \ref{Lemma1} that
$\mathnormal{\Phi} f\in N_{y_{\lambda}}$, that is, $\Phi
f(y_{\lambda})=0$ when $\lambda$ is large enough. Thus,
$\mathnormal{\Phi} f(y_{0})=0$ since $y_{\lambda}\rightarrow y_{0}$
and $\mathnormal{\Phi} f$ is continuous.

\par Let $\chi\in C(X)$ such that $\chi$ vanishes outside
$U_{\varphi(y_{0})}$ and $\chi(\varphi(y_{0}))=1$. Then, for any
$h$ in $C(X,E)$ we have $h=\chi h+(1-\chi)h$. Since $\chi h$
vanishes outside $U_{\varphi(y_{0})}$, by the above argument, we
can see that $\Phi(\chi h)(y_{0})=0$. Clearly,
$\mathnormal{\Phi}((1-\chi)h)$ vanishes at $y_{0}$ since
$(1-\chi)h\in M_{\varphi(y_{0})}$. Thus, $\mathnormal{\Phi}
h(y_{0})=\mathnormal{\Phi}(\chi
h)(y_{0})+\mathnormal{\Phi}((1-\chi)h)(y_{0})=0$ for any $h$ in
$C(X,E)$. This leads to a contradiction since $\mathnormal{\Phi}$
is surjective. So $\varphi$ is continuous, and thus a
homeomorphism from $Y$ onto $X$ satisfying $\mathnormal{\Phi}
M_{\varphi(y)} = N_{y}$ for each $y$ in $Y$.

Next, note that $\Ker \delta_{\varphi(y)}= \Ker
\delta_{y}\circ\mathnormal{\Phi}$, where $\delta_{y}$ is the Dirac
functional. Hence, there is a linear operator
$\mathnormal{\Pi}(y):E\rightarrow F$ such that
$\delta_{y}\circ\mathnormal{\Phi}=\mathnormal{\Pi}(y)\circ\delta_{\varphi(y)}$.
In other words,
$$
\mathnormal{\Phi} f(y)=\mathnormal{\Pi}(y)(f(\varphi(y))),\qquad \forall\,f\in C(X,
E),\,\forall\,y\in Y.
$$
See, e.g., \cite[p. 67]{8}.

It is a routine work to verify the other assertions in the statement
of this theorem. For the convenience of the readers, we give a
sketch of the rest of the proof. For $e$ in $E$, let
$\textbf{1}_{X}\otimes e\in C(X,E)$ be defined by
$(\textbf{1}_{X}\otimes e)(x)=e$  for each $x$ in $X$. Let $y$ in
$Y$ be fixed. If $e\neq0$, then
$\mathnormal{\Pi}(y)e=\mathnormal{\Pi}(y)((\textbf{1}_{X}\otimes
e)(\varphi(y)))=\mathnormal{\Phi}(\textbf{1}_{X}\otimes e)(y)\neq0$
since $\textbf{1}_{X}\otimes e$ is non-vanishing. Hence,
$\mathnormal{\Pi}(y)$ is one-to-one. On the other hand, for  $u$ in
$F$ we can find a function $f$ in $C(X,E)$ such that
$\mathnormal{\Phi} f=\textbf{1}_{Y}\otimes u$ by the surjectivity of
$\mathnormal{\Phi}$. Let $e=f(\varphi(y))$. Then
$\mathnormal{\Pi}(y)e=\mathnormal{\Pi}(y)(f(\varphi(y)))=\mathnormal{\Phi}
f(y)=u$. That is, $\mathnormal{\Pi}(y)$ is surjective. To see
$\mathnormal{\Pi}(y)$ is a Riesz isomorphism, let $e_{1}, e_{2}\in
E$. Then $\mathnormal{\Pi}(y)(e_{1}\vee
e_{2})=\mathnormal{\Phi}(\textbf{1}_{X}\otimes(e_{1}\vee
e_{2}))(y)=\mathnormal{\Phi}(\textbf{1}_{X}\otimes
e_{1})(y)\vee\mathnormal{\Phi}(\textbf{1}_{X}\otimes
e_{2})(y)=\mathnormal{\Pi}(y)e_{1}\vee\mathnormal{\Pi}(y)e_{2}$,
since $\mathnormal{\Phi}$ is a Riesz isomorphism.

\par Recall that every positive operator between Banach lattices is
continuous.  Let $e\in E$. Since
$\|\mathnormal{\Pi}(y)e\|=\|\mathnormal{\Phi}(\textbf{1}_{X}\otimes
e)(y)\|\leq\|\mathnormal{\Phi}(\textbf{1}_{X}\otimes
e)\|\leq\|\mathnormal{\Phi}\|\|e\|$, we have
$\|\mathnormal{\Pi}(y)\|\leq\|\mathnormal{\Phi}\|$ for all $y$ in $Y$.
On the other hand, for any $f$ in $C(X,E)$ and any $y$ in $Y$, we can
see $\|\mathnormal{\Phi}
f(y)\|=\|\mathnormal{\Pi}(y)(f(\varphi(y)))\|\leq\|\mathnormal{\Pi}(y)\|\|f\|$.
Consequently, $\|\mathnormal{\Phi}\|\leq \sup_{y\in Y} \|\mathnormal{\Pi}(y)\|$.

Finally, we prove that $\mathnormal{\Pi}:Y\rightarrow (\mathcal{L}(E, F), SOT)$ is
continuous. To this end, let $\{y_{\lambda}\}$ be a net such that
$y_{\lambda}\rightarrow y$ in $Y$. Then, for any $e$ in $E$,
$\|\mathnormal{\Pi}(y_{\lambda})e-\mathnormal{\Pi}(y)e\|=\|\mathnormal{\Phi}(\textbf{1}_{X}\otimes
e)(y_{\lambda})- \mathnormal{\Phi}(\textbf{1}_{X}\otimes
e)(y)\|\rightarrow 0$, since
$\mathnormal{\Phi}(\textbf{1}_{X}\otimes e)$ is continuous on
$Y$.
\end{proof}

\par In the above results, we have to assume that both $\mathnormal{\Phi}$ and
$\mathnormal{\Phi}^{-1}$  are  non-vanishing preserving. In the following
example, we can see that the inverse of a non-vanishing preserving Riesz
isomorphism is not necessarily non-vanishing preserving.

\begin{example}\label{Example1}
Let $X=\{1,2\}$ equipped with the discrete topology and
$E=\mathbb{R}$ with its usual ordering and norm, and let $Y=\{0\}$
and $F=\mathbb{R}^{2}$ with the pointwise ordering and the
 sup norm. Define $\mathnormal{\Phi}: C(X,E)\rightarrow C(Y,F)$ by
 $\mathnormal{\Phi}
f(0)=(f(1),\,f(2))$. Clearly, the Riesz isometric isomorphism
$\mathnormal{\Phi}$ is non-vanishing preserving, but its inverse
$\mathnormal{\Phi}^{-1}$ is not.
\end{example}

\par Let $E$, $F$ be both Banach lattices and Riesz
algebras,  Miao, Cao and Xiong \cite{14} recently proved that if
$F$ has no zero-divisor and there exists a Riesz algebraic
isomorphism $\mathnormal{\Phi}: C(X,E)\rightarrow C(Y,F)$ such
that $\mathnormal{\Phi}f$ is non-vanishing on $Y$ if $f$ is
non-vanishing on $X$, then $X$ is homeomorphic to $Y$, and $E$ is
Riesz algebraically isomorphic to $F$. In fact, from their proof
we can see that $\mathnormal{\Phi}f $ is non-vanishing on $Y$ if
and only if $f$ is non-vanishing on $X$, that is, both
$\mathnormal{\Phi}$ and $\mathnormal{\Phi}^{-1}$  are
non-vanishing preserving. Therefore, the result of Miao, Cao and
Xiong can be restated as follows.

\begin{corollary}[\cite{14}]\label{Corollary4}
Let $E$, $F$ be both Banach lattices and Riesz algebras. If $F$
has no zero-divisor and  $\mathnormal{\Phi}: C(X,E)\rightarrow
C(Y,F)$ is a Riesz algebraic isomorphism such that
$\mathnormal{\Phi}f $ is non-vanishing on $Y$ if $f$ is
non-vanishing on $X$, then $\mathnormal{\Phi}$ is a weighted
composition operator
$$
\mathnormal{\Phi} f(y)=\mathnormal{\Pi}(y)(f(\varphi(y))),\qquad \forall\,f\in
C(X, E),\,\forall\,y\in Y.
$$
Here $\varphi$ is a homeomorphism from
$Y$ onto $X$, and $\mathnormal{\Pi}(y)$ is a Riesz algebraic
isomorphism from $E$ onto $F$ for every $y$ in $Y$.
\end{corollary}

\par In Theorem \ref{Theorem2}, when $X$, $Y$ are compact Hausdorff
spaces and $E=F=\mathbb{R}$, the lattice hypothesis about
$\mathnormal{\Phi}$ can be dropped.

\begin{example}\label{eg:bi-invertibility-preserving}
Let $X,Y$ be compact Hausdorff spaces, and $C(X), C(Y)$ be the
Banach spaces of continuous real-valued functions defined on $X,
Y$, respectively. Assume $\mathnormal{\Phi}: C(X) \to C(Y)$ is a
linear map such that $\mathnormal{\Phi}f $ is non-vanishing on $Y$
if and only if $f$ is non-vanishing on $X$.

Note that
$(\mathnormal{\Phi}\mathbf{1}_{X})^{-1}\mathnormal{\Phi}$ is a
unital linear map preserving non-vanishing. Let $\lambda$ be in
the range of $f$. Then $f-\lambda \mathbf{1}_{X}$ is not
invertible, and thus neither is
$(\mathnormal{\Phi}\mathbf{1}_{X})^{-1}\mathnormal{\Phi}f -
\lambda \mathbf{1}_{Y}$.  It follows that $\lambda$ is in the
range of
$(\mathnormal{\Phi}\mathbf{1}_{X})^{-1}\mathnormal{\Phi}f$. The
converse also holds. Therefore, the range of
$(\mathnormal{\Phi}\mathbf{1}_{X})^{-1}\mathnormal{\Phi}f$
coincides with the range of $f$ for each $f$ in $C(X)$. In
particular, $(\mathnormal{\Phi}
\mathbf{1}_{X})^{-1}\mathnormal{\Phi}$ is a unital linear isometry
from $C(X)$ into $C(Y)$.  By the Holszty\'nski Theorem
\cite{Hol66}, there is a compact subset $Y_0$ of $Y$ and a
quotient map $\varphi : Y_0 \to X$ such that
$$
(\mathnormal{\Phi}\mathbf{1}_{X})^{-1}\mathnormal{\Phi}f\mid_{Y_0}
\ = f\circ\varphi, \qquad \forall f\in C(X).
$$

In case $\mathnormal{\Phi}$ is surjective, the classical
Banach-Stone Theorem ensures that $\varphi$ is a homeomorphism
from $Y=Y_0$ onto $X$. Moreover, if
$\mathnormal{\Phi}\mathbf{1}_{X}$ is strictly positive on $Y$,
then $\mathnormal{\Phi}$ is a Riesz isomorphism. However, when
$\mathnormal{\Phi}$ is not surjective the situation is a bit
uncontrollable. For example, consider $\mathnormal{\Phi} :
C[0,1]\to C([0,\frac{1}{2}]\cup[1,\frac{3}{2}])$ defined by
$$
\mathnormal{\Phi} f (y) = \left\{
            \begin{array}{ll}
              f(2y), & \text{if $0\leq y\leq 1/2$;} \\
              (2y-2)f(0) + (3-2y)f(1), & \text{if $1\leq y \leq \frac{3}{2}$.}
            \end{array}
          \right.
$$Clearly, the thus defined $\mathnormal{\Phi}$ is a
non-surjective linear isometry preserving non-vanishing in two
ways, but $[0,1]$ is not homeomorphic to
$[0,\frac{1}{2}]\cup[1,\frac{3}{2}]$.
\end{example}

Finally, we borrow an example from \cite{LW07} which shows that the
surjectivity cannot be guaranteed by many other properties we
usually consider.

\begin{example}\label{eg:beta}
Let $\omega$ and $\omega_{1}$ be the first infinite and the first
uncountable ordinal number, respectively. Let $[0,\omega_{1}]$ be
the compact Hausdorff space consisting of all ordinal numbers $x$
not greater than $\omega_{1}$ and equipped with the topology
generated by order intervals.  Note that every continuous function
$f$ in $C[0,\omega_{1}]$ is eventually constant.  More precisely,
there is a non-limit ordinal $x_f$ such that $\omega < x_f <
\omega_{1}$ and $f(x)=f(\omega_{1})$ for all $x \geq x_f$.

Define $\phi: [0,\omega_{1}]\to [0,\omega_{1}]$ by setting
$$
\mathnormal{\phi}(0)=\omega_{1}, \quad \mathnormal{\phi}(n)=n-1
\text{ for all $n = 1, 2, \ldots$, \quad and }
\mathnormal{\phi}(x)=x \text{ for all } x\geq \omega.
$$
Let $\mathnormal{\Phi} : C[0,\omega_{1}]\to C[0,\omega_{1}]$ be
the \emph{non-surjective} composition operator defined by
$\mathnormal{\Phi}f = f\circ\mathnormal{\phi}$. It is plain that
$\mathnormal{\Phi}$ is an isometric unital algebraic and lattice
isomorphism from $C[0,\omega_{1}]$ onto its range. In fact, one
can see in \cite[Example 3.3]{LW07} that the map
$\mathnormal{\Phi}$  is a non-surjective linear $n$-local
automorphism of $C[0,\omega_{1}]$, where $n=1,2,\ldots,\omega$,
i.e., the action of $\mathnormal{\Phi}$ on any set $S$ of
cardinality not greater than $n$ agrees with an automorphism
$\mathnormal{\Phi}_S$.

\end{example}
\section*{Acknowledgment}
The authors would like to thank  the referee for his comments
which have improved this paper.

\end{document}